\documentclass[11pt
]{amsart}

\usepackage{hyperref}
\hypersetup{bookmarksdepth=3}

\usepackage{geometry}
\usepackage{amsmath,amssymb}
\usepackage{
mathrsfs}
\usepackage{esint}

\usepackage{tensor}

\usepackage{amssymb,amsmath,amsthm}

\newtheorem{theorem}{Theorem}

\theoremstyle{remark}\newtheorem{remark}[theorem]{Remark}

\usepackage{graphicx}

\begin{document}

\title{Continuity of weighted estimates for sublinear operators}
\author{Michael Papadimitrakis and Nikolaos Pattakos}

\date{}

\begin{abstract}
{In this note we prove that if a sublinear operator $T$ satisfies a certain weighted estimate in the $L^{p}(w)$ space for all $w\in A_{p}$, $1<p<+\infty$, then 

$$\lim_{d_{*}(w,w_{0})\to0}\|T\|_{L^{p}(w)\rightarrow L^{p}(w)}=\|T\|_{L^{p}(w_{0})\rightarrow L^{p}(w_{0})},$$
where $d_{*}$ is the metric defined in \cite{NV} and $w_{0}$ is a fixed $A_{p}$ weight . }
\end{abstract}

\maketitle

\begin{section}{introduction and notation}
We are going to work with positive $L^{1}_{loc}(\mathbb R^n)$ functions $w$ (called weights), that satisfy the following condition for some $1<p<+\infty$,

$$[w]_{A_{p}}:=\sup_{Q}\Big(\frac1{|Q|}\int_{Q}w(x)dx\Big)\Big(\frac1{|Q|}\int_{Q}w(x)^{-\frac1{p-1}}dx\Big)^{p-1}<+\infty.$$
The number $[w]_{A_{p}}$ is called the $A_{p}$ characteristic of the weight $w$ and we say that $w\in A_{p}$. The supremum is taken over all cubes $Q$ of $\mathbb R^n$. 

In \cite{NV} the authors defined a metric $d_{*}$ in the set of $A_{p}$ weights. For two weights $u,v\in A_{p}$ we define

$$d_{*}(u,v):=\|\log u-\log v\|_{*},$$
where for a function $f$ in $L^{1}_{loc}(\mathbb R^n)$ we define the $BMO(\mathbb R^n)$ norm (modulo constants) as

$$\|f\|_{*}:=\sup_{Q}\frac1{|Q|}\int_{Q}|f(x)-f_{Q}|dx.$$
The notation $f_{Q}$ is used to denote the average value of the function $f$ over the cube $Q$ (we will also use the notation $<f>_{Q}$). In addition, the authors proved that if a linear operator $T$ satisfies the weighted estimate

$$\|T\|_{L^{p}(w)\rightarrow L^{p}(w)}\leq F([w]_{A_{p}}),$$
for all $w\in A_{p}$, where $F$ is a positive increasing function, then for any fixed weight $w_{0}\in A_{p}$ we have

$$\lim_{d_{*}(w,w_{0})\to0}\|T\|_{L^{p}(w)\rightarrow L^{p}(w)}=\|T\|_{L^{p}(w_{0})\rightarrow L^{p}(w_{0})},$$
which means that the operator norm of $T$ on the $L^{p}(w)$ space is a continuous function of the weight $w$ with respect to the $d_{*}$ metric. In this note we are going to extend this result for sublinear operators $T$. Namely, we have the

\begin{theorem}
\label{main}
Suppose that for some $1<p<+\infty$, a sublinear operator $T$ satisfies the inequality

$$\|T\|_{L^{p}(w)\rightarrow L^{p}(w)}\leq F([w]_{A_{p}}),$$
for all $w\in A_{p}$, where $F$ is a positive increasing function. Fix an $A_{p}$ weight $w_{0}$. Then

$$\lim_{d_{*}(w,w_{0})\to0}\|T\|_{L^{p}(w)\rightarrow L^{p}(w)}=\|T\|_{L^{p}(w_{0})\rightarrow L^{p}(w_{0})}.$$
\end{theorem}

Let us mention that the method used in \cite{NV} can not be used for sublinear operators. The argument there does not work for them.

\begin{remark}
In \cite{B} Buckley showed that the Hardy-Littlewood maximal operator defined as

$$Mf(x)=\sup_{x\in Q}\frac1{|Q|}\int_{Q}|f(y)|dy,$$
where the supremum is taken over all cubes $Q$ in $\mathbb R^n$ that contain the point $x$, satisfies the estimate

$$\|M\|_{L^{p}(w)\rightarrow L^{p}(w)}\leq c[w]_{A_{p}}^{\frac1{p-1}},$$
for $1<p<+\infty$, and all weights $w\in A_{p}$, where the constant $c>0$ is independent of the weight $w$. This means that the assumptions of Theorem \ref{main} hold for $M$. 
\end{remark}   

We present the proof of the Theorem in the next section.

\end{section}

\begin{section}{proof of theorem \ref{main}}
The main tool for the proof is the inequality (proved in \cite{NV}) 

\begin{equation}
\label{in1}
\|T\|_{L^{p}(u)\rightarrow L^{p}(u)}\leq\|T\|_{L^{p}(v)\rightarrow L^{p}(v)}(1+c_{[v]_{A_{p}}}d_{*}(u,v)),
\end{equation}
that holds for all $A_{p}$ weights $u,v\in A_{p}$ that are sufficiently close in the $d_{*}$ metric, and for sublinear operators $T$ that satisfy the assumptions of our Theorem. The positive constant $c_{[v]_{A_{p}}}$ that appears in the inequality depends on the dimension $n$, $p$, the function $F$ and the $A_{p}$ characteristic of the weight $v$. Since the quantities $n, p, F$ are fixed we only write the subscript $c_{[v]_{A_{p}}}$ to emphasize this dependence on the characteristic.

\begin{proof}
We apply inequality (\ref{in1}) with $u=w$ and $v=w_{0}$ to obtain

$$\|T\|_{L^{p}(w)\rightarrow L^{p}(w)}\leq\|T\|_{L^{p}(w_{0})\rightarrow L^{p}(w_{0})}(1+c_{[w_{0}]_{A_{p}}}d_{*}(w,w_{0})).$$
By letting $d_{*}(w,w_{0})$ go to $0$ we get

$$\limsup_{d_{*}(w,w_{0})\to0}\|T\|_{L^{p}(w)\rightarrow L^{p}(w)}\leq\|T\|_{L^{p}(w_{0})\rightarrow L^{p}(w_{0})}.$$
Now it suffices to prove the inequality

$$\|T\|_{L^{p}(w_{0})\rightarrow L^{p}(w_{0})}\leq\liminf_{d_{*}(w,w_{0})\to0}\|T\|_{L^{p}(w)\rightarrow L^{p}(w)},$$
in order to finish the proof. For this reason we use inequality (\ref{in1}) with $u=w_{0}$ and $v=w$

$$\|T\|_{L^{p}(w_{0})\rightarrow L^{p}(w_{0})}\leq\|T\|_{L^{p}(w)\rightarrow L^{p}(w)}(1+c_{[w]_{A_{p}}}d_{*}(w,w_{0})).$$ 
At this point if we know that the constant $c_{[w]_{A_{p}}}$ remains bounded as the distance $d_{*}(w,w_{0})$ goes to $0$ we are done. 

For this reason we assume that $d_{*}(w,w_{0})=\delta$ is very close to $0$. Then the function $\frac{w}{w_{0}}$ is an $A_{p}$ weight with $A_{p}$ characteristic very close to $1$ (see \cite{GCRF}). How close depends only on $\delta$, not on $w$. Thus, if $R$ is large enough, the weight $(\frac{w}{w_{0}})^{R}\in A_{p}$, with $A_{p}$ characteristic independent of $w$ (again see \cite{GCRF}). Note that from the classical $A_{p}$ theory, for sufficiently small $\epsilon>0$, we have $w_{0}^{1+\epsilon}\in A_{p}$. Choose the numbers $R, \epsilon$ such that we have the relation $\frac{1}{R}+\frac{1}{1+\epsilon}=1$, i.e. such that $R$ and $R'=1+\epsilon$ are conjugate numbers. Then, by applying H\"older's inequality twice we have the following

\begin{eqnarray*}
<w>_{Q}<w^{-\frac1{p-1}}>_{Q}^{p-1}&=&\Big<\frac{w}{w_{0}}w_{0}\Big>_{Q}\Big<\Big(\frac{w}{w_{0}}\Big)^{-\frac1{p-1}}w_{0}^{-\frac1{p-1}}\Big>_{Q}^{p-1}\\
&\leq&\Big<\Big(\frac{w}{w_{0}}\Big)^{R}\Big>_{Q}^{\frac1{R}}\Big<w_{0}^{R'}\Big>_{Q}^{\frac1{R'}}\Big<\Big(\frac{w}{w_{0}}\Big)^{-\frac1{p-1}\cdot R}\Big>_{Q}^{\frac{p-1}{R}}\Big<w_{0}^{-\frac1{p-1}\cdot R'}\Big>_{Q}^{\frac{p-1}{R'}}\\
&=&\Big<\Big(\frac{w}{w_{0}}\Big)^{R}\Big>_{Q}^{\frac1{R}}\Big<\Big(\Big(\frac{w}{w_{0}}\Big)^{R}\Big)^{-\frac1{p-1}}\Big>_{Q}^{\frac{p-1}{R}}\Big<w_{0}^{R'}\Big>_{Q}^{\frac1{R'}}\Big<(w_{0}^{R'})^{-\frac1{p-1}}\Big>_{Q}^{\frac{p-1}{R'}}\\
&\leq&\Big[\Big(\frac{w}{w_{0}}\Big)^{R}\Big]_{A_{p}}^{\frac1{R}}[w_{0}^{1+\epsilon}]_{A_{p}}^{\frac1{R'}}\leq C,
\end{eqnarray*}
where $C$ is a constant independent of the weight $w$. Therefore, $[w]_{A_{p}}\leq C$. 

The last step is to remember how we obtain the constant $c_{[w]_{A_{p}}}$ that appears in inequality (\ref{in1}). The authors in \cite{NV} used the Riesz-Thorin interpolation theorem with change in measure and then expressed one of the terms that appears in their calculations as a Taylor series. The constant $c_{[w]_{A_{p}}}$ appears at exactly this point and it is not difficult to see that it depends continuously on $[w]_{A_{p}}$. Since this characteristic is bounded for $w$ close to $w_{0}$ in the $d_{*}$ metric we have that $c_{[w]_{A_{p}}}$ is bounded as well. This completes the proof.

\end{proof} 

A consequence of the proof is the following remark.

\begin{remark}
Fix a weight $w_{0}\in A_{p}$ and a positive number $\delta$ sufficiently small. There is a positive constant $C$ that depends on $[w_{0}]_{A_{p}}$ and $\delta$ such that for all weights $w$ with $d_{*}(w,w_{0})<\delta$ we have $[w]_{A_{p}}\leq C$. In addition, from the inequality (see the proof of Theorem \ref{main})

\begin{equation}
\label{in2}
[w]_{A_{p}}\leq\Big[\Big(\frac{w}{w_{0}}\Big)^{R}\Big]_{A_{p}}^{\frac1{R}}[w_{0}^{1+\epsilon}]_{A_{p}}^{\frac1{R'}},
\end{equation}
and Lebesgue dominated convergence theorem (by letting $R\to+\infty$ and remembering that the $A_{p}$ constant of the weight $(\frac{w}{w_{0}})^{R}$ is independent of $R$) we obtain

$$\limsup_{d_{*}(w,w_{0})\to 0}[w]_{A_{p}}\leq[w_{0}]_{A_{p}}.$$ 
In order to get the remaining inequality

$$[w_{0}]_{A_{p}}\leq\liminf_{d_{*}(w,w_{0})\to 0}[w]_{A_{p}},$$
we rewrite (\ref{in2}) as

$$[w_{0}]_{A_{p}}\leq\Big[\Big(\frac{w_{0}}{w}\Big)^{R}\Big]_{A_{p}}^{\frac1{R}}[w^{1+\epsilon}]_{A_{p}}^{\frac1{R'}},$$
and we proceed in the same way as before. In this case the number $\epsilon$ depends on $[w]_{A_{p}}$. But we already know that for $w$ close to $w_{0}$ in the $d_{*}$ metric the $A_{p}$ characteristic of $w$ is bounded from above. This means that we are allowed to choose the same number $\epsilon$ for all weights $w$ that are sufficiently close to $w_{0}$ and we are done. Therefore, the $A_{p}$ characteristic of a weight $w\in A_{p}$ is a continuous function of the weight with respect to the $d_{*}$ metric.

\end{remark}

\end{section} 

\textbf{Acknowledgments}: The authors would like to thank professor Alexander Volberg from Michigan State University in East Lansing for useful discussions.

M. Papadimitrakis, Department of Mathematics, University of Crete, Knossos Ave., 71409 Iraklion, Greece;
papadim@math.uoc.gr;\newline

N. Pattakos, Department of Mathematics, Michigan State University, East Lansing, Michigan 48824, USA;
pattakos@msu.edu;

\end{document}